\documentclass[12pt]{article}

\usepackage{amsmath,amssymb,amsthm}
\usepackage {graphicx}

\newtheorem{theorem}{Theorem}[section]
\newtheorem{lemma}[theorem]{Lemma}
\newtheorem{remark}[theorem]{Remark}

\newtheorem{problem}[theorem]{Problem}
\newtheorem {corollary}[theorem]{Corollary}
\def \RR{\mathbb R}
\def \NN{\mathbb N}
\def\({\left(}  \def\){\right)}
\def\[{\left[}  \def\]{\right]}

\def \beq {\begin {equation}}
\def \eeq {\end{equation}}
\def \OL {\overline}

\def \W {\widetilde}
\def\NN {\mathbb N}

\begin {document}
\begin {center}
\Large{
Optimal Antipodal Configuration of $2d$ Points \\ 
on a Sphere in $\mathbb R^d$ for Covering}

\bigskip

Sergiy Borodachov

\bigskip

\small {\it Department of Mathematics, Towson University, 8000 York Rd., Towson, MD, 21252}

\end {center}

\begin {abstract}
We show that among antipodal $2d$-point configurations on the sphere $S^{d-1}$ in $\RR^d$, the set of vertices of a regular cross-polytope inscribed in $S^{d-1}$ uniquely solves the best-covering problem (this is new for $d\geq 5$) and the maximal polarization problem for potentials given by a function of the distance squared with a positive and convex second derivative ($d\geq 3$).
\end {abstract}


\medskip


\large{
\section {Introduction}

Let $S^{d-1}=\{(x_1,\ldots,,x_{d})\in \RR^d: x_1^2+\ldots+x_{d}^2=1\}$, $d\geq 2$, be the unit sphere in the Euclidean space $\RR^d$. For a given point configuration $\omega_N:=\{{\bf x}_1,\ldots,{\bf x}_N\}\subset S^{d-1}$, 
let 
$$
\rho(\omega_N,S^{d-1}):=\max\limits_{{\bf x}\in S^{d-1}}\min\limits_{i=\OL{1,N}}\left|{\bf x}-{\bf x}_i\right|
$$
denote its {\it mesh norm} (relative to the sphere). The classical optimal covering problem on the sphere is to minimize the mesh-norm.
\begin {problem}\label {covering'}
{\rm
For a given $N\in \NN$, find the quantity
\begin {equation}\label {p1}
\rho_N(S^{d-1}):=\inf\limits_{\omega_N\subset S^{d-1}}\rho(\omega_N,S^{d-1})
\end {equation}
and optimal-covering $N$-point configurations on $S^{d-1}$; i.e., configurations $\omega^\ast_N\subset S^{d-1}$ attaining the infimum on the right-hand side of \eqref {p1}. }
\end {problem}

Solution to Problem \ref {covering'} is known on $S^1$ for every $N$ (equally spaced points), on $S^2$ for $N=1-8$, $10$, $12$, and $14$, see the works by Fejes-T\'oth \cite {Tot1953,Tot1964,Tot1969}, Sch\"utte \cite {Sch1955}, and Wimmer \cite {Wim2017}, on $S^3$ for $N=1-6$ and $8$, see the works by Galiev \cite {Gal1996}, B\"or\"oczky and Wintsche \cite {BorWin2003}, and Dalla, Larman, Mani-Levitska, and Zong \cite {DalLarMan2000}, and on $S^{d-1}$, $d>4$, for $1\leq N\leq d+2$, see \cite {Gal1996,BorWin2003}. For more information on optimal covering, see, among others, books by Fejes-T\'oth \cite {Tot1953,Tot1964}, Rogers \cite {Rog2008}, and B\"or\"oczky \cite{Bor2004}. 

The main goal of this paper is a further study of the optimal-covering property of the set of vertices of a {\it regular cross-polytope} inscribed in $S^{d-1}$, which is any $2d$-point configuration of the form $\OL\omega_{2d}:=\{\pm {\bf a}_1,\ldots,\pm{\bf a}_d\}$, where $\{{\bf a}_1,\ldots,{\bf a}_d\}$ is an orthonormal basis in $\RR^d$. The optimal covering property of $\OL\omega_{2d}$ is known on $S^{d-1}$ for $d=2$ (a basic result), $d=3$ \cite{Tot1953,Tot1964}, and $d=4$ \cite {DalLarMan2000}. It is a well-known open question for $d\geq 5$. We resolve it here for antipodal configurations. 

The proof in \cite {DalLarMan2000} for $d=4$ uses the area argument together with the Euler-Poincar\'e formula
and the fact that (cf. \cite {Bor1987}) the largest ($d-1$)-dimensional area of a spherical simplex inscribed in a given spherical cap of angular radius less than $\pi/2$ is that of a regular simplex. If the number $\tau$ of ($d-1$)-dimensional simplices obtained after triangulating the facets of the convex hull of a given $2d$-point configuration is less than or equal to $2^d$, the area argument combined with the above mentioned fact from \cite {Bor1987} complete the proof. If $\tau>2^d$ then the Euler-Poincar\'e formula is used to provide an efficient estimate for the number of one-dimensional edges of the convex hull. Since the number of terms in the Euler-Poincar\'e formula grows with $d$, such an estimate is not possible when $d\geq 5$ and the argument fails. However, if one makes an additional assumption that the configurations are antipodal, then it can be shown that $\tau=2^d$, and the area argument can be used together with the above mentioned fact from \cite {Bor1987} to complete the proof, see Theorem~\ref {covering}. 

The paper is structured as follows. We state our main result, Theorem~\ref {covering}, in Section \ref {two} and its consequences in Section~\ref {twoa}. Section \ref {three} contains the proof of Theorem \ref {covering}. In Section \ref {foura} we state and prove Theorem~\ref {mstiff} on the minimum value of the potential of a regular cross-polytope, which we use in Section \ref {four} to prove consequences of Theorem \ref {covering}. The known proof of Lemma \ref {p} is given in Appendix (Section \ref {five}) for completeness.

\section {Main result}\label {two}

Recall that a point configuration on $S^{d-1}$ is called {\it antipodal} if together with a point ${\bf x}$ it contains $-{\bf x}$. We say that a given point set in $\RR^d$ is {\it in general position} if it is not contained in any hyperplane (or, equivalently, in any $(d-1)$-dimensional affine subspace). 
It will be convenient for us to recast the optimal-covering problem in terms of dot products.
\begin {remark}
{\rm 
A point configuration $\omega_N\subset S^{d-1}$ is a solution to Problem \ref {covering'} if and only if $\omega_N$ maximizes the quantity
$$
\eta(\omega_N,S^{d-1}):=\min\limits_{{\bf x}\in S^{d-1}}\max\limits_{i=\OL{1,N}}{\bf x}\cdot {\bf x}_i
$$
over all $N$-point configurations $\omega_N=\{{\bf x}_1,\ldots,{\bf x}_N\}$ on $S^{d-1}$.
}
\end {remark}

We allow configurations, where points may coincide. However, if an antipodal configuration of $2d$ points on $S^{d-1}$ is in general position, its points are pairwise distinct.

The main result of this paper is the following. 
\begin {theorem}\label {covering}
Let $d\geq 2$ and $\omega_{2d}$ be an antipodal configuration of $2d$ points on $S^{d-1}$. Then
\begin {equation}\label {3}
\eta(\omega_{2d},S^{d-1})\leq \frac {1}{\sqrt{d}}.
\end {equation}
Equality in \eqref {3} holds if and only if $\omega_{2d}$ is the set of vertices of a regular cross-polytope inscribed in~$S^{d-1}$.
\end {theorem}
For $d=3$ and $4$, Theorem \ref {covering} is a special case of theorems proved in \cite {Tot1953,Tot1964} and in \cite {DalLarMan2000}, respectively. For $d=2$, it is a basic result.

\section {A consequence of Theorem \ref {covering} for polarization}\label {twoa}

The best-covering problem is a limiting case of the polarization problem stated next. Let $f:[0,4]\to(-\infty,\infty]$ be a function finite and continuous on $(0.4]$ such that $f(0)=\lim\limits_{t\to 0^+}f(t)$. We will call $f$ a {\it potential function}. For a given point configuration $\omega_N=\{{\bf x}_1,\ldots,{\bf x}_N\}\subset S^{d-1}$, 
denote
$$
P_f(\omega_N,S^{d-1}):=\min\limits_{{\bf x}\in S^{d-1}}\sum\limits_{i=1}^{N}f\(\left|{\bf x}-{\bf x}_i\right|^2\)
$$
and let
\begin {equation}\label {2}
\mathcal P_f(S^{d-1},N):=\sup\limits_{\omega_N\subset S^{d-1}}P_f(\omega_N,S^{d-1}).
\end {equation}
The (max-min) $N$-point polarization problem on the sphere is stated in the following way.
\begin{problem}\label {P2}
{\rm
Find quantity
\eqref {2}
and $N$-point configurations $\omega_N^\ast$ on $S^{d-1}$ that attain the supremum on the right-hand side of~\eqref {2}.
}
\end {problem}
In the case of a sphere, solution to Problem \ref {P2} is known on the unit circle $S^1$ for every $N\geq 1$ and $f(t)=t^{-s/2}$, $s>0$, and $f(t)=-t^{-s/2}$, $-1\leq s<0$, which correspond to the Riesz potential, as well as for $f(t)=\frac {1}{2}\ln \frac {1}{t}$, which corresponds to the logarithmic potential, see the works by Stolarsky, Ambrus, Nikolov, Rafailov, Ball, Erd\'elyi, Saff, Hardin, and Kendall \cite {Sto1975circle,Amb2009,NikRaf2011,AmbBalErd2013,ErdSaf2013,HarKenSaf2013}. The solution is the set of vertices of a regular $N$-gon inscribed in $S^1$.
In fact, paper \cite {HarKenSaf2013} established this for arbitrary interactions given by a decreasing and convex function of the geodesic distance on $S^1$. Optimal $N$-point configurations for polarization on $S^1$ were also characterized by Bosuwan and Ruengrot \cite {BosRue2017} for the Riesz potential with $s=-2,-4,\ldots,2-2N$, $N\geq 2$. 

On the sphere $S^{d-1}$, $d\geq 3$, the solution to Problem \ref {P2} is known for $1\leq N\leq d$ (basic result) and $N=d+1$, see the works by Su \cite {Su2014} and the author \cite {Bor2022} showing the optimality of a regular simplex. Also, for $N=2d$, the optimality of a regular cross-polytope $\OL\omega_{2d}$ was shown among centered configurations by Boyvalenkov, Dragnev, Hardin, Saff, and Stoyanova 
\cite {BoyDraHarSafSto600cell}. A configuration $\omega_{2d}\subset S^{d-1}$ is called {\it centered} if there is a point ${\bf y}\in S^{d-1}$ such that $-\frac {1}{\sqrt{d}}\leq {\bf y}\cdot {\bf x}_i\leq \frac {1}{\sqrt{d}}$, $i=1,\ldots,2d$.

Other settings of polarization problem were studied in \cite {Sto1975circle,Mon1988,ErdSaf2013,BosRue2017,AmbNie2019,BetFauSte,Borminmax,BoyDraHarSafSto600cell}. 
More extensive reviews on polarization (including the continuous version and asymptotics) can be found, for example, in book \cite {BorHarSaf2019}.

Polarization problem can also be recast in terms of dot products.
Let $g(t):=f(2-2t)$, $t\in [-1,1]$. Then $f$ is a potential function if and only if $g:[-1,1]\to (-\infty,\infty]$ is a function finite and continuous on $[-1,1)$ such that $g(1)=\lim\limits_{t\to 1^-}g(t)$. For every configuration $\omega_N=\{{\bf x}_1,\ldots,{\bf x}_N\}\subset S^{d-1}$, we have
$$
P^g(\omega_N,S^{d-1}):=\min\limits_{{\bf x}\in S^{d-1}}\sum\limits_{i=1}^{N}g({\bf x}\cdot{\bf x}_i)=\min\limits_{{\bf x}\in S^{d-1}}\sum\limits_{i=1}^{N}f\(\left|{\bf x}-{\bf x}_i\right|^2\)=P_f(\omega_N,S^{d-1}).
$$
Thus, a point configuration $\omega_N$ maximizes the quantity $P_f(\omega_N,S^{d-1})$ if and only if it maximizes the quantity $P^g(\omega_N,S^{d-1})$.

Inequality \eqref {3} in Theorem \ref {covering} is equivalent to the following statement.
\begin {corollary}\label {centered'}
Any antipodal $2d$-point configuration on $S^{d-1}$, $d\geq 2$, is centered.
\end {corollary}

Paper \cite {BoyDraHarSafSto600cell} proved universal bounds for polarization and used them to show that regular cross-polytope $\OL\omega_{2d}$ is optimal for polarization among all centered configurations on $S^{d-1}$. This result together with Corollary \ref {centered'} imply the following (the equality in \eqref {P} below follows from Theorem \ref {mstiff}).

\begin {corollary}\label {polarization}
Suppose $g:[-1,1]\to (-\infty,\infty]$ is a function continuous on $[-1,1)$ with $g(1)=\lim\limits_{t\to 1^-}g(t)$ and differentiable on $(-1,1)$ such that $g''$ is non-negative and convex on $(-1,1)$. Suppose also that $\omega_{2d}\subset S^{d-1}$, $d\geq 2$, is any antipodal $2d$-point configuration. Then 
\begin {equation}\label {P}
P^g(\omega_{2d},S^{d-1})\leq P^g(\OL\omega_{2d},S^{d-1})=d\(g\(\frac {1}{\sqrt{d}}\)+g\(-\frac {1}{\sqrt{d}}\)\).
\end {equation}
\end {corollary}
Here we provide an elementary proof of Corollary \ref {polarization} that avoids the use of universal bounds for polarization and is valid for any centered configuration. Concerning uniqueness of the optimal configuration, we have the following.
\begin {corollary}\label {uniqueness}
Under assumptions of Corollary \ref {polarization}, if $g''$ is strictly positive on $(-1,1)$, then equality holds throughout \eqref {P} if and only if $\omega_{2d}$  is the set of vertices of a regular cross-polytope inscribed in $S^{d-1}$.
\end {corollary}

We remark that a similar uniqueness question for centered configurations remains open. Observe that we do not require the potential function $g$ to be monotone in Corollaries \ref {polarization} and \ref {uniqueness}.

\section {Proof of Theorem \ref {covering}}\label {three}

We need the following auxiliary statement proved by B\"or\"oczky in \cite {Bor1987} (see also \cite [Lemma 6.7.2]{Bor2004}). Let $r:\RR^d\setminus \{{\bf 0}\}\to S^{d-1}$, $r({\bf x})={\bf x}/\left|{\bf x}\right|$, denote the radial projection onto $S^{d-1}$. We will call a {\it $(d-1)$-simplex} a simplex in $\RR^d$ with $d$ vertices having dimension $(d-1)$.


%

\begin {lemma}\label {p}
Let $H$ be a hyperplane in $\RR^d$, $d\geq 2$, with equation $x_d=a$, where $0<a<1$. Let $Y$ be a non-degenerate $(d-1)$-simplex inscribed in $H\cap S^{d-1}$. Then the $(d-1)$-dimensional volume of the radial projection $r(Y)$ of $Y$ onto $S^{d-1}$ is the largest if and only if $Y$ is a regular $(d-1)$-simplex inscribed in $H\cap S^{d-1}$.
\end {lemma}
For completeness, we provide the proof of Lemma \ref {p} in the Appendix.

Let $D(\omega_{2d})$ denote the convex hull of a point configuration $\omega_{2d}\subset S^{d-1}$. Points from $\omega_{2d}$ will also be called vertices.
Recall that a hyperplane ${\bf x}\cdot {\bf a}=\alpha$ in $\RR^d$ is called a {\it supporting hyperplane for a convex body $B$} if for every ${\bf z}\in B$, we have ${\bf z}\cdot {\bf a}\leq \alpha$, while for some ${\bf y}\in B$, we have ${\bf y}\cdot {\bf a}=\alpha$.

\begin {lemma}\label {q}
Let $\omega_{2d}\subset S^{d-1}$, $d\geq 2$, be an antipodal configuration of $2d$ points in general position. Then the interior of $D(\omega_{2d})$ contains the origin, and the boundary of $D(\omega_{2d})$ is the union of $2^d$ $(d-1)$-simplices whose vertices are in $\omega_{2d}$, pairwise intersections have $(d-1)$-dimensional measure~$0$, and the hyperplane containing each simplex does not pass through the origin. 
\end {lemma}
\begin {proof}
Pick an arbitrary subset of $d$ linearly independent vectors from $\omega_{2d}$ and denote it by $\{{\bf y}_1,\ldots,{\bf y}_d\}$. Such a subset exists, since $\omega_{2d}$ is in general position.
Let $\Pi:=\{-1,1\}^d$. Then any of the $2^d$ sets $T_{\boldsymbol\sigma}:=\{\sigma_1{\bf y}_1,\ldots,\sigma_d{\bf y}_d\}$, where $\boldsymbol\sigma=(\sigma_1,\ldots,\sigma_d)\in \Pi$, is a linearly independent subset of $\omega_{2d}$ and, hence, is a set of vertices of a $(d-1)$-simplex, which we will denote by $F_{\boldsymbol\sigma}$. For every $\boldsymbol\sigma\in \Pi$, the set $T_{\boldsymbol\sigma}$ is contained in a unique hyperplane $H_{\boldsymbol\sigma}$. Then the set $T_{-\boldsymbol\sigma}=-T_{\boldsymbol\sigma}$ is contained in the hyperplane $-H_{\boldsymbol\sigma}$ with $-H_{\boldsymbol\sigma}\neq H_{\boldsymbol\sigma}$, since $\omega_{2d}=T_{\boldsymbol\sigma}\cup T_{-\boldsymbol\sigma}$ is in general position. Then $D(\omega_{2d})$ is contained in the closed subset of $\RR^d$ bounded by $H_{\boldsymbol\sigma}$ and $-H_{\boldsymbol\sigma}$. Then the simplex $F_{\boldsymbol\sigma}$ is contained in the boundary $\partial D(\omega_{2d})$ of $D(\omega_{2d})$ for every $\boldsymbol\sigma\in \Pi$.

Since $\omega_{2d}$ is in general position, the origin is in the interior of $D(\omega_{2d})$ as a point on the line segment joining two points in relative interiors of opposite facets. 

Let now ${\bf x}$ be any point in $\partial D(\omega_{2d})$. Then ${\bf x}\neq {\bf 0}$ and there is a unique set of numbers $\alpha_1,\ldots,\alpha_d$ such that ${\bf x}=\alpha_1{\bf y}_1+\ldots+\alpha_d{\bf y}_d$. For some $\boldsymbol\sigma=(\sigma_1,\ldots,\sigma_d)\in \Pi$, we have ${\bf x}=\beta_1\sigma_1{\bf y}_1+\ldots+\beta_d\sigma_d{\bf y}_d$, where $\beta_i:=\left|\alpha_i\right|\geq 0$, $i=1,\ldots,d$. Let $\beta:=\sum_{i=1}^{d}\beta_i$. Then $\beta>0$ and $(1/\beta){\bf x}\in F_{\boldsymbol\sigma}\subset \partial D(\omega_{2d})$. If it were that $\beta>1$, then, since $(1/\beta){\bf x}\in H_{\boldsymbol\sigma}$, the origin and ${\bf x}$ would lie in different half-spaces relative to $H_{\boldsymbol\sigma}$, which is a supporting hyperplane for $D(\omega_{2d})$. Since the origin is in $D(\omega_{2d})$, we have ${\bf x}\notin D(\omega_{2d})$; that is, ${\bf x}\notin \partial D(\omega_{2d})$. If it were that $\beta<1$, then ${\bf x}$ would be in the relative interior of the line segment joining $(1/\beta){\bf x}$ and ${\bf 0}$. Since $(1/\beta){\bf x}$ is in $D(\omega_{2d})$ and ${\bf 0}$ is in its interior, then ${\bf x}$ is also in the interior of $D(\omega_{2d})$; that is, ${\bf x}\notin \partial D(\omega_{2d})$. This contradiction shows that $\beta=1$. Then ${\bf x}\in F_{\boldsymbol\sigma}$.

Thus, $\partial D(\omega_{2d})=\bigcup\limits_{\boldsymbol\sigma\in \Pi}F_{\boldsymbol\sigma}$. Assume that the intersection of two simplices $F_{{\boldsymbol\sigma}}$ and $F_{{\boldsymbol\sigma}'}$, with ${\boldsymbol\sigma}=(\sigma_1,\ldots,\sigma_d)\neq {\boldsymbol\sigma}'=(\sigma_1',\ldots,\sigma_d')$ is non-empty. For every point ${\bf z}$ in both simplices, there are numbers $\beta_1,\ldots,\beta_d,\beta_1',\ldots,\beta_d'\geq 0$ such that $\sum_{i=1}^{d}\beta_i=\sum_{i=1}^{d}\beta_i'=1$ and 
$$
{\bf z}=\beta_1\sigma_1{\bf y}_1+\ldots+\beta_d\sigma_d{\bf y}_d=\beta_1'\sigma_1'{\bf y}_1+\ldots+\beta_d'\sigma_d'{\bf y}_d.
$$
Since vectors ${\bf y}_1,\ldots,{\bf y}_d$ are linearly independent, we have $\beta_i\sigma_i=\beta_i'\sigma_i'$, $i=1,\ldots,d$. Since $\sigma_j\neq \sigma_j'$ for some $j$, we have $\beta_j=-\beta_j'$. In view of the non-negativity, we have $\beta_j=\beta_j'=0$. Therefore, ${\bf z}$ belongs to a lower-dimensional face of $F_{{\boldsymbol\sigma}}$. Thus, ${\rm Vol}_{d-1}(F_{{\boldsymbol\sigma}}\cap F_{{\boldsymbol\sigma'}})=0$. Furthermore, the hyperplane $H_{\boldsymbol\sigma}$ containing $F_{\boldsymbol\sigma}$ does not contain the origin, since if it did, vectors ${\bf y}_1,\ldots,{\bf y}_d$ would be linearly dependent.
\end {proof}

Recall that an intersection $U$ of a convex polytope $P$ with its supporting hyperplane, such that ${\rm dim}\ \! U=d-1$, is called a {\it facet} of $P$. The polytope $D(\omega_{2d})$ in Lemma \ref {q} has $2^d$ facets. 

\begin {proof}[Proof of Theorem \ref {covering}]
We first verify that $\eta(\OL\omega_{2d},S^{d-1})=\frac {1}{\sqrt{d}}$. Without loss of generality, we can assume that $\OL\omega_{2d}=\W\omega_{2d}:=\{\pm {\bf e}_1,\ldots,\pm {\bf e}_d\}$, where ${\bf e}_1,\ldots,{\bf e}_d$ are standard basis vectors in $\RR^d$. Then for every vector ${\bf x}=(x_1,\ldots,x_d)\in S^{d-1}$, we have $\max\limits_{i=\OL{1,d}}{\bf x}\cdot (\pm {\bf e}_i)=\max\limits_{i=\OL{1,d}}\left|x_i\right|\geq \frac {1}{\sqrt {d}}$, since $\sum_{i=1}^{d}\left|x_i\right|^2=1$, with equality occuring whenever each coordinate of ${\bf x}$ is $\pm\frac {1}{\sqrt{d}}$.

Assume to the contrary that there is an antipodal configuration $\omega_{2d}\subset S^{d-1}$ such that $\eta(\omega_{2d},S^{d-1})> 1/\sqrt{d}$. Then $D(\omega_{2d})$ contains ${\bf 0}$ in its interior (if it did not, we would have $\eta(\omega_{2d},S^{d-1})\leq 0<1/\sqrt{d}$). Furthermore, $\omega_{2d}$ is in general position and, hence, its points are pairwise distinct. Polytope $D(\omega_{2d})$ contains the sphere $S$ of radius $1/\sqrt{d}$ centered at the origin. If it didn't, then any point ${\bf a}\in S\setminus D(\omega_{2d})$ would be strictly separated from $D(\omega_{2d})$ by some hyperplane $L=\{{\bf x} : {\bf x}\cdot {\bf v}=c\}$, where we can take $\left|{\bf v}\right|=1$ and ${\bf v}\cdot {\bf a}>c$. Then for every ${\bf x}_i\in \omega_{2d}$, we would have ${\bf x}_i\cdot {\bf v}<c<{\bf a}\cdot {\bf v}\leq \left|{\bf a}\right|\left|{\bf v}\right|= 1/\sqrt{d}$ which would contradict the contrary assumption. 

By Lemma \ref {q}, boundary of $D(\omega_{2d})$ is the union of $(d-1)$-simplices (which we denote by $F_1,\ldots,F_{2^d}$) with vertices in $\omega_{2d}$. Their pairwise intersections have $(d-1)$-dimensional volume $0$, and there are $a_i>0$ and ${\bf z}_i\in S^{d-1}$, $i=1,\ldots,2^d$, such that the hyperplane $H_i:=\{{\bf x} : {\bf x}\cdot {\bf z}_i=a_i\}$ contains $F_i$. Denote by $k$ an index such that the radial projection $r(F_k)$ onto $S^{d-1}$ of the simplex $F_k$ has $(d-1)$-dimensional volume at least $1/2^{d}$ of the $(d-1)$-dimensional volume of $S^{d-1}$. If $V$ denotes a regular $(d-1)$-simplex inscribed in $H_k\cap S^{d-1}$, then in view of Lemma \ref {p}, we have 
$$
2^{-d}{\rm Vol}_{d-1}(S^{d-1})\leq {\rm Vol}_{d-1}(r(F_k))\leq {\rm Vol}_{d-1}(r(V)). 
$$
Since $H_k$ does not contain the origin (by Lemma \ref {q}), there are no antipodal pairs among the vertices of $F_k$. Then the remaining $d$ points from $\omega_{2d}$ are contained in the hyperplane ${\bf x}\cdot {\bf z}_k=-a_k$. We have $a_k=\max\limits_{i=\OL {1,2d}}{\bf z}_k\cdot{\bf x}_i\geq \eta(\omega_{2d},S^{d-1})>1/\sqrt{d}$. The $(d-1)$-dimensional volume (denoted by $\nu$) of the radial projection $r(W)$ of a regular simplex $W$ incribed in $M\cap S^{d-1}$, where $M=\{{\bf x} : {\bf x}\cdot{\bf z}_k=1/\sqrt{d}\}$, will be strictly larger than ${\rm Vol}_{d-1}(r(V))$. This is because $a_k>1/\sqrt{d}$ and the radius of the intersection of a hyperplane perpendicular to ${\bf z}_k$ with $S^{d-1}$ decreases as the hyperplane moves further away from the origin. Since $r(W)$ is the same as the radial projection of any facet of $D(\OL\omega_{2d})$ and there are exactly $2^d$ facets, we have $$
2^{-d}{\rm Vol}_{d-1}(S^{d-1})\leq {\rm Vol}_{d-1}(r(V))<{\rm Vol}_{d-1}(r(W))=2^{-d}{\rm Vol}_{d-1}(S^{d-1}). 
$$
This contradiction proves \eqref {3}.

To complete the proof of Theorem \ref {covering}, assume that equality holds in $\eqref {3}$. Then the origin is in the interior of $D(\omega_{2d})$. Assume to the contrary that there is a facet $Q$ of $D(\omega_{2d})$ such that the hyperplane containing it is at a distance strictly greater than $1/\sqrt{d}$ from the origin. Then by Lemma \ref {p}, the radial projection of $Q$ onto $S^{d-1}$ has $(d-1)$-dimensional volume strictly less than $\nu$. Equality in \eqref {3} implies that every other facet $Y$ of $D(\omega_{2d})$ is contained in a hyperplane whose distance to the origin is at least $1/\sqrt{d}$. By Lemma \ref {p}, ${\rm Vol}_{d-1}(r(Y))\leq \nu$ and, by Lemma \ref {q}, $D(\omega_{2d})$ has $2^d$ facets. Then the radial projection onto $S^{d-1}$ of the whole boundary of $D(\omega_{2d})$ has $(d-1)$-dimensional volume strictly less than $2^d\nu={\rm Vol}_{d-1}(S^{d-1})$. This contradiction shows that the hyperplane containing each facet of $D(\omega_{2d})$ is at a distance exactly $1/\sqrt{d}$ from the origin. Each facet of $D(\omega_{2d})$ must be a regular simplex. If some facet $U$ were not, by Lemma \ref {p}, we would have ${\rm Vol}_{d-1}(r(U))<\nu$ while for any other facet $J$, we would have ${\rm Vol}_{d-1}(r(J))\leq\nu$ leading to a similar contradiction. Let ${\bf a}_1,\ldots,{\bf a}_d$ be the vertices of one of the facets of $D(\omega_{2d})$. Then $\omega_{2d}=\{\pm {\bf a}_1,\ldots,\pm {\bf a}_d\}$. It is not difficult to verify that ${\bf a}_i\cdot {\bf a}_j=0$, $i\neq j$. Then $\{{\bf a}_1,\ldots,{\bf a}_d\}$ is an orthonormal basis in $\RR^d$ and $\omega_{2d}$ is the set of vertices of a regular cross-polytope inscribed in~$S^{d-1}$.
\end {proof}

\section {Minimum of the potential of a regular cross-polytope}\label {foura}

In order to show the optimality of $\OL\omega_{2d}$ for Problem \ref {P2}, one will need to know the quantity $P_f(\OL\omega_{2d},S^{d-1})$ by locating the absolute minima of the potential of $\OL\omega_{2d}$ on the sphere $S^{d-1}$.
This was done earlier for the vertices of a regular $N$-gon inscribed in $S^1$ and for the vertices of a regular simplex, cross-polytope, and cube inscribed in $S^{d-1}$ for Riesz potential functions ($s\neq 0$) and their horizontal translations, see the works by Stolarsky, Nikolov, and Rafailov \cite {Sto1975circle,Sto1975,NikRaf2011,NikRaf2013}. For general potentials, this has been recently done in \cite {Bor2022} for vertices of a regular simplex. Below, we extend one of the results of \cite {Sto1975,NikRaf2013} by finding the absolute minima of the potential of $\OL\omega_{2d}$ for potential functions $g$ that have a convex second derivative. Our proof is different from the one in \cite {Sto1975,NikRaf2013} and uses polynomial interpolation and convexity of $g''$. Without loss of generality, we can assume in this section that $\OL\omega_{2d}=\W\omega_{2d}=\{\pm {\bf e}_1,\ldots,\pm {\bf e}_{d}\}$, where ${\bf e}_1,\ldots,{\bf e}_d$ is the standard basis in $\RR^d$.
\begin {theorem}\label {mstiff}
Let $d\geq 2$ and $g:[-1,1]\to(-\infty,\infty]$ be a function continuous on $[-1,1)$ and differentiable on $(-1,1)$ such that $g(1)=\lim\limits_{t\to 1^-}g(t)$ and $g''$ is convex on $(-1,1)$. Then the potential 
$$
p^g(\W\omega_{2d},{\bf x}):=\sum\limits_{{\bf y}\in \W\omega_{2d}}g({\bf x}\cdot {\bf y})
$$
achieves its absolute minimum over $S^{d-1}$ at any point of $S^{d-1}$ whose every coordinate is $1/\sqrt{d}$ or $-1/\sqrt{d}$ (these points are vertices of the cube dual to $\W\omega_{2d}$). Furthermore,
$$
P^g(\W\omega_{2d},S^{d-1})=d\(g\(\frac {1}{\sqrt{d}}\)+g\(-\frac {1}{\sqrt{d}}\)\).
$$
\end {theorem}
We also remark that in the upcoming paper \cite {Borstiff}, we obtain the locations of absolute minima of the potential of any configuration on $S^{d-1}$ which is a tight spherical design of an even strength or a ($2m-1$)-design contained in the union of $m$ parallel hyperplanes. Such is, for example, $\OL\omega_{2d}$ for $m=2$.

\begin {proof}[Proof of Theorem \ref {mstiff}]
Let $h(t):=g(t)+g(-t)$. Then $h$ is even and $h''$ is convex on $(-1,1)$. Let $p$ be the Hermite interpolating polynomial for $h$ at points $t_1=-1/\sqrt{d}$ and $t_2=1/\sqrt{d}$. Then the even polynomial $(p(t)+p(-t))/2$ is also Hermite for $h$ at $t_1$ and $t_2$. By uniqueness, $p$ must be even. Since $p$ has degree at most $3$, it has the form $p(t)=at^2+b$. Furthermore, $\frac {a}{d}+b=p(1/\sqrt{d})=h\(1/\sqrt{d}\)$.

We also have, $h(t)\geq p(t)$, $t\in [-1,1]$. Indeed, assume to the contrary that $v(t):=h(t)-p(t)$ is negative for some $t=t_0\in (-1,1)$. Note that $v(t_1)=v(t_2)=v'(t_1)=v'(t_2)=0$. Then the Mean value theorem and the Rolle's theorem imply that there are points $-1<\tau_1<\tau_2<\tau_3<1$ such that $v''(\tau_1)<0$ and $v''(\tau_2)=v''(\tau_3)=0$ (if $t_0<t_1$) or $v''(\tau_1)<0$, $v''(\tau_2)>0$, and $v''(\tau_3)<0$ (if $t_1<t_0<t_2$) or $v''(\tau_1)=v''(\tau_2)=0$ and $v''(\tau_3)<0$ (if $t_0>t_2$). None of these cases is possible, since $h''$ (and, hence, $v''$) is convex on $(-1,1)$. Thus, $v$ is non-negative on $(-1,1)$; that is, $h(t)\geq p(t)$, $t\in (-1,1)$. We extend this inequality to the endpoints by passing to the limit. 

Let ${\bf x}=(x_1,\ldots,x_d)$ be an arbitrary point on $S^{d-1}$. Then
\begin{equation*}
\begin {split}
p^g(\W\omega_{2d},{\bf x})&=\sum\limits_{i=1}^{d}\(g(x_i)+g(-x_i)\)=\sum\limits_{i=1}^{d}h(x_i)\geq \sum\limits_{i=1}^{d}p(x_i)=\sum\limits_{i=1}^{d}(ax_i^2+b)\\
&=a+d\cdot b=d\cdot h\(\frac {1}{\sqrt{d}}\)=d\(g\(\frac{1}{\sqrt{d}}\)+g\(-\frac {1}{\sqrt{d}}\)\)=p^g(\W\omega_{2d},{\bf x}^\ast),
\end {split}
\end {equation*}
where ${\bf x}^\ast$ is any point on $S^{d-1}$ whose every coordinate is $1/\sqrt{d}$ or $-1/\sqrt{d}$.
\end {proof}

\section {Proof of Corollaries \ref {polarization} and \ref {uniqueness}}\label {four}

\begin {proof}[Proof of Corollary \ref {polarization}]We choose an arbitrary centered configuration $\omega_{2d}=\{{\bf y}_1,\ldots,{\bf y}_{2d}\}\subset S^{d-1}$
and let ${\bf x}^\ast\in S^{d-1}$ be a point such that 
\begin {equation}\label {centered}
-\frac {1}{\sqrt{d}}\leq {\bf x}^\ast \cdot{\bf y}_i\leq \frac {1}{\sqrt{d}}, \ \ \ i=1,\ldots,2d. 
\end {equation}
In fact, $-\frac {1}{\sqrt{d}}\leq \pm {\bf x}^\ast \cdot{\bf y}_i\leq \frac {1}{\sqrt{d}}$ for all $i$. Let $-t_1\leq \ldots\leq-t_{2d}\leq t_{2d}\leq \ldots\leq t_1$ denote all the dot products that ${\bf x}^\ast$ and $-{\bf x}^\ast$ form with points ${\bf y}_i$. Since $g''$ is non-negative on $(-1,1)$, $g'$ is non-decreasing on $(-1,1)$, and, hence, the function $h(t)=g(t)+g(-t)$ is non-decreasing on $[0,1)$.  Since $t_1,\ldots,t_{2d}\in [0,1/\sqrt{d}]$, we have
\begin {equation}\label {q1}
\begin {split}
P^g&(\omega_{2d},S^{d-1})\leq \frac {1}{2}\sum\limits_{i=1}^{2d}g({\bf x}^\ast\cdot {\bf y}_i)+\frac {1}{2}\sum\limits_{i=1}^{2d}g(-{\bf x}^\ast\cdot {\bf y}_i)=\frac {1}{2}\sum\limits_{i=1}^{2d}\(g(t_i)+g(-t_i)\)\\
&=\frac {1}{2}\sum\limits_{i=1}^{2d}h(t_i)\leq d\cdot h\(\frac {1}{\sqrt{d}}\)=d\(g\(\frac {1}{\sqrt{d}}\)+g\(-\frac {1}{\sqrt{d}}\)\)=P^g(\OL\omega_{2d},S^{d-1}),
\end {split}
\end {equation}
where the last equality in \eqref {q1} holds in view of Theorem \ref {mstiff}. This proves \eqref {P} for any centered configuration $\omega_{2d}$. In view of Corollary \ref {centered'}, we have Corollary~\ref {polarization}.
\end {proof}

\begin {proof}[Proof of Corollary \ref {uniqueness}] Assume that equality holds throughout \eqref {P} for a given antipodal configuration $\omega_{2d}\subset S^{d-1}$. Since $\omega_{2d}$ is antipodal, ${\bf x}^\ast$ can be chosen in the beginning of the proof of Corollary \ref {polarization} with the additional property that 
$\max\limits_{i=\OL {1,2d}}{\bf x}^\ast\cdot {\bf y}_i=\eta(\omega_{2d},S^{d-1})$ (inequalities \eqref {centered} will then hold in view of Theorem~\ref {covering}). Since equality now holds throughout \eqref {q1}, we have $\sum_{i=1}^{2d}h(t_i)=2d\cdot h(1/\sqrt{d})$. Since $0\leq t_i\leq 1/\sqrt{d}$, $i=1,\ldots,2d$ and $h$ is strictly increasing on $[0,1)$, we have $t_i=1/\sqrt{d}$, $i=1,\ldots,2d$.
Then ${\bf x}^\ast$ forms only dot products $1/\sqrt{d}$ and $-1/\sqrt{d}$ with points of $\omega_{2d}$. Thus, $\eta(\omega_{2d},S^{d-1})=\max\limits_{i=\OL {1,2d}}{\bf x}^\ast\cdot {\bf y}_i=1/\sqrt{d}$. By the uniqueness part of Theorem \ref {covering},  configuration $\omega_{2d}$ is the set of vertices of a regular cross-polytope inscribed in~$S^{d-1}$.
\end {proof}


\section {Appendix. Proof of Lemma \ref {p}.}\label {five}

The assertion of Lemma \ref {p} is trivial for $d=2$. Therefore, we assume that $d\geq 3$.
Let $Y$ be a $(d-1)$-simplex inscribed in $H\cap S^{d-1}$ whose radial projection $r(Y)$ onto $S^{d-1}$ has the largest $(d-1)$-dimensional volume. Let ${\bf y}_1,\ldots,{\bf y}_d$ be the vertices of $Y$. Let $\Omega(Y)$ be the intersection of the closed unit ball $B^d$ centered at the origin with the cone, denoted by ${\rm cone}\{{\bf y}_1,\ldots,{\bf y}_d\}$, constructed as the convex hull of the rays starting at ${\bf 0}$ and passing through each vertex of $Y$. Then ${\rm Vol}_{d-1}(r(Y))=d\cdot {\rm Vol}_d(\Omega(Y))$. 

Assume to the contrary that $Y$ is not a regular simplex. Then there is a vertex of $Y$ such that some two edges stemming out of it have non-equal lengths. Without loss of generality, we can assume that $\left|{\bf y}_3-{\bf y}_1\right|\neq \left|{\bf y}_3-{\bf y}_2\right|$. Let $L$ be the hyperplane that is the perpendicular bisector for the line segment with endpoints ${\bf y}_1$ and ${\bf y}_2$. Observe that ${\bf 0}\in L$. Denote by ${\bf y}_i'$ the orthogonal projection of the point ${\bf y}_i$ onto $L$, $i=3,\ldots,d$, and let $Y'$ be the convex hull of $\{{\bf y}_1,{\bf y}_2,{\bf y}_3',\ldots,{\bf y}_d'\}$. Since ${\bf y}_1,{\bf y}_2,{\bf y}_3',\ldots,{\bf y}_d'\in H\cap B^d$, we have $Y'\subset H\cap B^d$ and $r(Y')\subset C:=\{(x_1,\ldots,x_d)\in S^{d-1} : x_d\geq a\}$. Also, $r(Y')={\rm cone}\{{\bf y}_1,{\bf y}_2,{\bf y}_3',\ldots,{\bf y}_d'\}\cap S^{d-1}$. 

Recall that a Steiner symmetrization of a set $A\subset \RR^d$ relative to a hyperplane $L$ is the set 
$$
{\rm St}(A):=\bigcup\limits_{\ell\in \mathcal Q}\(\frac {1}{2}(A\cap \ell)+\frac {1}{2}(\W A\cap \ell)\), 
$$
where $\mathcal Q$ is the set of all lines passing through points of $A$ that are perpendicular to hyperplane $L$ and $\W A$ is the reflection of $A$ with respect to $L$. The Steiner symmetrization of ${\rm cone}\{{\bf y}_1,\ldots,{\bf y}_d\}$ is contained in ${\rm cone}\{{\bf y}_1,{\bf y}_2,{\bf y}_3',\ldots,{\bf y}_d'\}$. Since $L$ passes through ${\bf 0}$, whenever $A,D\subset \RR^d$ are such that ${\rm St}(A)\subset D$, we have ${\rm St}(A\cap B^d)\subset D\cap B^d$. Therefore, ${\rm St}(\Omega(Y))\subset \Omega(Y')$, where $\Omega(Y'):={\rm cone}\{{\bf y}_1,{\bf y}_2,{\bf y}_3',\ldots,{\bf y}_d'\}\cap B^d$. Since Steiner symmetrization preserves the volume, we have 
\begin {equation*}
\begin {split}
{\rm Vol}_{d-1}(r(Y))&={d}{{\rm Vol}_d(\Omega(Y))}={d}{{\rm Vol}_d[{\rm St}(\Omega(Y))]}\\&\leq {d}{{\rm Vol}_d(\Omega(Y'))}={\rm Vol}_{d-1}(r(Y')).
\end {split}
\end {equation*}
Since the $(d-1)$-dimensional volume of $r(Y)$ is positive, so is the one of $r(Y')$ and, hence of $Y'$. Then $Y'$ is a $(d-1)$-simplex.
Since ${\bf y}_3\notin L$, its projection ${\bf y}'_3$ is an interior point of $H\cap B^d$. Then there is a point ${\bf y}_3''\in H\cap B^d$ into which we can move ${\bf y}'_3$ so that the new simplex $Y''$ with vertices ${\bf y}_1,{\bf y}_2,{\bf y}_3'',{\bf y}_4'\ldots,{\bf y}_d'$ contains $Y'$ together with some open set that lies outside of $Y'$. Then we will have ${\rm Vol}_{d-1}(r(Y'))<{\rm Vol}_{d-1}(r(Y''))$ with $Y''\subset H\cap B^d$. The simplex $Y''$ is contained in some $(d-1)$-simplex $Y'''$ inscribed in $H\cap S^{d-1}$. Then ${\rm Vol}_{d-1}(r(Y))<{\rm Vol}_{d-1}(r(Y'''))$. This contradicts the assumption that $r(Y)$ has the largest volume over all $(d-1)$-simplices $Y$ inscribed in $H\cap S^{d-1}$.

Thus, if ${\rm Vol}_{d-1}(r(Y))$ is maximal, then $Y$ must be a regular simplex. At the same time, for all regular $(d-1)$-simplices $V$ inscribed in $H\cap S^{d-1}$, the volumes ${\rm Vol}_{d-1}(r(V))$ are the same. Then ${\rm Vol}_{d-1}(r(V))$ is maximal.

}

\end {document}